\documentclass[11pt]{article}

\usepackage{amssymb,amsmath,amsthm,enumitem,enumitem,mathabx,caption}

\newtheorem{thm}{Theorem}[section]
\newtheorem{prp}[thm]{Proposition}

\newtheorem{crl}[thm]{Corollary}

\theoremstyle{definition}
\newtheorem{dfn}[thm]{Definition}

\theoremstyle{remark}

\def\lra{\longrightarrow}

\def\BE#1{\begin{equation}\label{#1}}
\def\EE{\end{equation}}
\def\lr#1{\langle#1\rangle}
\def\flr#1{\left\lfloor{#1}\right\rfloor}
\def\blr#1{\big\langle#1\big\rangle}

\def\wt#1{\widetilde{#1}}
\def\ov#1{\overline{#1}}
\def\eref#1{(\ref{#1})}
\def\tn#1{\textnormal{#1}}
\def\sf#1{\textsf{#1}}

\def\wh#1{\widehat{#1}}

\def\fc{\mathfrak c}
\def\ff{\mathfrak f}
\def\fj{\mathfrak j}

\def\om{\omega}
\def\si{\sigma}
\def\vph{\varphi}
\def\ze{\zeta}

\def\Ga{\Gamma}
\def\La{\Lambda}
\def\Si{\Sigma}

\def\C{\mathbb C}
\def\cJ{\mathcal J}
\def\cM{\mathcal M}
\def\fM{\mathfrak M}
\def\P{\mathbb P}
\def\R{\mathbb R}
\def\Q{\mathbb Q}
\def\cT{\mathcal T}
\def\Z{\mathbb Z}

\def\bu{\bullet}
\def\dbar{\bar\partial}
\def\i{\infty}

\def\tnd{\tn{d}}
\def\rdet{\wh{\tn{det}}}
\def\ev{\tn{ev}}
\def\E{\tn{E}}
\def\GW{\tn{GW}}
\def\Hom{\tn{Hom}}
\def\id{\tn{id}}
\def\pt{\tn{pt}}

\def\top{\tn{top}}

\begin{document}

\title{Real Orientations, Real Gromov-Witten Theory,\\
and Real Enumerative Geometry}
\author{Penka Georgieva\thanks{Supported by ERC grant STEIN-259118} $~$and 
Aleksey Zinger\thanks{Partially supported by NSF grant DMS 1500875 and MPIM}}
\date{\today}
\maketitle

\begin{abstract}
\noindent
The present note overviews our recent construction of real Gromov-Witten theory
in arbitrary genera for many real symplectic manifolds, 
including the odd-dimensional projective spaces and the renowned quintic threefold,
its properties, and its connections with real enumerative geometry.
Our construction  introduces the principle of orienting
the determinant of a differential operator relative to a suitable base operator
and a real setting analogue of the (relative) spin structure of open Gromov-Witten theory.
Orienting the relative determinant, which  in the now-standard cases
is canonically equivalent  to orienting the usual determinant,
is naturally related to the topology of vector bundles in the relevant category.
This principle and its applications allow us to endow the uncompactified moduli 
spaces of real maps from symmetric surfaces of all topological types 
with natural orientations and to verify that they extend 
across the codimension-one boundaries of these spaces,
thus implementing a far-reaching proposal from C.-C.~Liu's thesis. 
\end{abstract}

\section{Real maps}
\label{RealSt_maps}

\noindent
The study of curves in projective varieties has been central to algebraic geometry
since the nineteenth century.
It was reinvigorated through its introduction
into symplectic topology in~\cite{Gr} and 
now plays prominent roles in symplectic topology and string theory as~well.
The foundations of (complex) \sf{Gromov-Witten} (or \sf{GW-}) \sf{theory}, 
i.e.~of counts of $J$-holomorphic curves in symplectic manifolds, 
were established in the 1990s and have been spectacularly applied ever~since.
On the other hand, the progress in establishing the foundations of \sf{real GW-theory},
i.e.~of counts of $J$-holomorphic curves in symplectic manifolds preserved by 
anti-symplectic involutions, has been much slower:
it did not exist in positive genera until~\cite{RealGWsI}.\\

\noindent
A \sf{real symplectic manifold} is a triple $(X,\om,\phi)$ consisting 
of a symplectic manifold~$(X,\om)$ and an anti-symplectic involution~$\phi$.
For such a triple, we denote by ${\mathcal J}_{\om}^{\phi}$ the space of $\om$-compatible
 almost complex structures~$J$ on~$X$ such that $\phi^*J\!=\!-J$.
The fixed locus~$X^{\phi}$ of~$\phi$ is then a Lagrangian submanifold of~$(X,\om)$
which is totally real with respect to any $J\!\in\!\cJ_{\om}^{\phi}$.
The basic example of a real Kahler manifold $(X,\om,\phi,J)$ is
the complex projective space~$\P^{n-1}$ with the Fubini-Study symplectic form,
the coordinate conjugation
$$\tau_n\!: \P^{n-1}\lra \P^{n-1}, \qquad 
\tau_n\big([Z_1,\ldots,Z_n]\big)=\big[\ov{Z}_1,\ldots,\ov{Z}_n\big],$$
and the standard complex structure.
Another important example is a real \textsf{quintic threefold}~$X_5$,
i.e.~a smooth hypersurface in~$\P^4$ cut out by a real equation;
it plays a prominent role in the interactions with string theory and algebraic geometry.\\

\noindent
A \sf{symmetric surface} $(\Si,\si)$ 
is a connected oriented, possibly nodal, surface~$\Si$ with an orientation-reversing involution~$\si$.
A \sf{symmetric Riemann surface} $(\Si,\si,\fj)$ 
is a symmetric surface $(\Si,\si)$ with an almost complex structure~$\fj$ on~$\Si$
such that $\si^*\fj\!=\!-\fj$.
A continuous map 
$$u\!:(\Si,\si)\lra(X,\phi)$$ 
is called \sf{real} if $u\!\circ\!\si\!=\!\phi\!\circ\!u$.
Such a map is said to be \sf{of degree} $B\!\in\!H_2(X;\Z)$ if $u_*[\Si]\!=\!B$.
Real maps~$u$ from~$(\Si,\si,\fj)$ and $u'$ from~$(\Si',\si',\fj')$ to~$(X,\phi)$
are \sf{equivalent} if there exists a biholomorphic map $h\!:\Si\!\lra\!\Si'$
such that $u\!=\!u'\!\circ\!h$, $j\!=\!h^*j'$, and $h\!\circ\!\si\!=\!\si'\!\circ\!h$.\\

\noindent
There are $\flr{\frac{3g+4}{2}}$ topological types of {\it smooth} symmetric surfaces; 
see \cite[Corollary~1.1]{Nat}.
As described in \cite[Section~3]{Melissa}, 
there are four types of nodes a 
one-nodal symmetric surfaces  $(\Si,x_{12},\si)$ may~have
\begin{enumerate}[leftmargin=.3in]

\item[(E)]\label{E_it} 
$x_{12}$ is an isolated real node, i.e.~$x_{12}$ is an isolated point
of the fixed locus $\Si^{\si}\!\subset\!\Si$;

\item[(H)]  $x_{12}$ is a non-isolated real node and

\begin{enumerate}[label=(H\arabic*),leftmargin=.15in] 

\item the topological component $\Si_{12}^{\si}$ of $\Si^{\si}$ containing $x_{12}$
is algebraically irreducible (the normalization~$\wt\Si_{12}^{\wt\si}$
of~$\Si_{12}^{\si}$ is connected);

\item the topological component $\Si_{12}^{\si}$ of $\Si^{\si}$ containing $x_{12}$
is algebraically reducible, but $\Si$ is algebraically irreducible
(the normalization~$\wt\Si_{12}^{\wt\si}$ of~$\Si_{12}$ is disconnected,
but the normalization $\wt\Si$ of~$\Si$ is connected);

\item  $\Si$ is algebraically reducible (the normalization $\wt\Si$ of~$\Si$ is disconnected).

\end{enumerate}
\end{enumerate} 
In the genus~0 case, the degenerations~(E) and~(H3) are known as 
the  \sf{codimension~1 sphere bubbling} and \sf{disk bubbling}, respectively;
the degenerations~(H1) and~(H2) cannot occur in the genus~0 case.
The one-nodal symmetric surfaces can be smoothed out in one-parameter families
to symmetric surfaces, typically  of different involution types for the two directions
of smoothings (the smoothings of~(H3) are always of the same type though).

\section{Moduli spaces of real maps}
\label{RealMod_maps}

\noindent 
Let $(X,\om,\phi)$ be a real symplectic manifold, $g,l\!\in\!\Z^{\ge0}$, 
$B\!\in\!H_2(X;\Z)$, and $J\!\in\!\cJ_{\om}^{\phi}$.
For a smooth symmetric surface $(\Si,\si)$, we denote by 
\BE{notgluedspace_e}\fM_{g,l}(X,B;J)^{\phi,\si}\subset \ov\fM_{g,l}(X,B;J)^{\phi,\si}\EE
the uncompactified moduli space of degree~$B$ real $J$-holomorphic maps from 
$(\Si,\si)$ to $(X,\phi)$ with $l$~conjugate pairs of marked points
and its stable map compactification.
A (virtually) codimension~1 stratum of $\ov\fM_{g,l}(X,B;J)^{\phi,\si}$
consists of maps from one-nodal symmetric surfaces.
By the existence of precisely two directions of smoothings of one-nodal symmetric surfaces,
each such stratum is either a hypersurface in $\ov\fM_{g,l}(X,B;J)^{\phi,\si}$ 
or a boundary of the spaces $\ov\fM_{g,l}(X,B;J)^{\phi,\si}$ 
for precisely two topological types of orientation-reversing involutions~$\si$ on~$\Si$. 
Thus, the union of real moduli spaces
\BE{gluedspace_e}
\ov\fM_{g,l}(X,B;J)^{\phi}=\bigcup_{\si}\ov\fM_{g,l}(X,B;J)^{\phi,\si}\EE
over all topological types of orientation-reversing involutions~$\si$ on~$\Si$ 
forms a space without boundary.
If $g\!+\!l\!\ge\!2$, there is a natural forgetful morphism
\BE{ffdfn_e}\ff\!: \ov\fM_{g,l}(X,B;J)^{\phi}\lra \R\ov\cM_{g,l}\equiv \ov\fM_{g,l}(\pt,0)^{\id}\EE
to the Deligne-Mumford moduli space of marked real curves.
\\

\noindent
The uncompactified moduli spaces in complex GW-theory have canonical 
orientations; see \cite[Section~3.2]{MS}.
As the ``boundary" strata of the moduli spaces in complex GW-theory have real codimension
of at least~2, this orientation automatically extends over the entire moduli space.
The two main difficulties in developing real GW-theory is
the potential non-orientability of $\fM_{g,l}(X,B;J)^{\phi,\si}$
and the fact that its virtual boundary strata have real codimension~1.
The origins of real GW-theory go back to~\cite{Melissa},
where the spaces~\eref{gluedspace_e} are topologized by
adapting the description of Gromov's topology in~\cite{LT}
via versal families of deformations of abstract complex curves to the real setting.
This  demonstrates that the codimension~1 boundaries of the spaces 
in~\eref{notgluedspace_e}  form hypersurfaces inside 
the full moduli space~\eref{gluedspace_e} and thus reduces the problem of constructing 
a real GW-theory to showing~that
\begin{enumerate}[label=(\Alph*),leftmargin=*]

\item\label{orient_it} the uncompactified moduli spaces $\fM_{g,l}(X,B;J)^{\phi,\si}$
are orientable for all types of orientation-reversing involutions~$\si$ on 
a smooth genus~$g$ symmetric surface, and 

\item\label{bnd_it} an orientation on
$$\fM_{g,l}(X,B;J)^{\phi}=\bigcup_{\si}\fM_{g,l}(X,B;J)^{\phi,\si}$$
extends across the (virtually) codimension~1 strata of $\ov\fM_{g,l}(X,B;J)^{\phi}$.\\ 

\end{enumerate} 

\noindent
Invariant counts of real curves were first constructed in~\cite{Wel4,Wel6}
following a different approach.
They were defined only in genus~0, for real symplectic 4- and 6-folds,
 and under certain topological conditions 
ruling out maps from type~(E) nodal symmetric surfaces
(thus, only (H3) degenerations can occur).
As they concerned only primary constraints, they did not give rise 
to a fully fledged real GW-{\it theory}.
The relevant moduli spaces in the settings of~\cite{Wel4,Wel6} are in fact not orientable,
and the invariance of the defined counts is checked by following  the paths
of curves induced by
paths between two generic almost complex  structures and two generic collections of
constraints.
In the interpretation of~\cite{Wel4,Wel6} in~\cite{Cho,Sol}, 
this invariance corresponds to the relevant moduli spaces being orientable
outside of (virtual) hypersurfaces not crossed by the paths of stable maps
induced by paths between two generic almost  complex structures and 
two generic collections of constraints.\\

\noindent
A fully fledged real GW-theory in genus~0 was finally set up in~\cite{Ge2}
following the original approach in~\cite{Melissa} and 
establishing~\ref{orient_it} and~\ref{bnd_it}  under certain topological conditions
on $(X,\om,\phi)$ and the map degree~$B$.
The topological conditions in~\cite{Ge2} ruling out maps from type~(E) 
nodal symmetric surfaces were later removed in~\cite{Teh}.
The genus~0 real GW-theory of \cite{Ge2,Teh} is used in~\cite{RealEnum}
to establish a real analogue of the WDVV relation of~\cite{KM,RT}
by pulling back a relation on $\R\ov\cM_{0,3}$ by the forgetful morphism~\eref{ffdfn_e}.\\

\noindent
The perspective on orienting the relevant moduli spaces taken in \cite{Cho,Sol,Ge2,Teh} is 
heavily influenced by the approach in the open GW-theory going back to
the late 1990s and the initial version of~\cite{FOOO}.
It works well in genus~0 because a splitting of a smooth genus~0 symmetric surface~$(\Si,\si)$
into two bordered surfaces interchanged by~$\si$ is unique up to homotopy and
the one-nodal transitions between the (two) different involution types in genus~0 
preserve such a splitting.
The former is no longer the case for most smooth symmetric surfaces of genus $g\!\ge\!2$;
the latter is not the case for most transitions in genus $g\!\ge\!1$.
This means that understanding the orientability of the moduli spaces of maps from
bordered (half-) surfaces is not sufficient for understanding the orientability 
of the moduli spaces of real maps and a new perspective to this problem is needed.
Such a perspective,  which is intrinsically real, rather than a ``doubled open",
is introduced in the construction of all-genera real GW-theory 
in~\cite{RealGWsI} and is summarized in Sections~\ref{RealOrient_sec} and~\ref{MainPf_sec} 
of this note.

\section{Real orientations}
\label{RealOrient_sec}

\noindent
Let $(X,\phi)$ be a topological space with an involution.
A \sf{conjugation} on a complex vector bundle $V\!\lra\!X$ 
\sf{lifting} an involution~$\phi$ is a vector bundle homomorphism 
$\vph\!:V\!\lra\!V$ covering~$\phi$ (or equivalently 
a vector bundle homomorphism  $\vph\!:V\!\lra\!\phi^*V$ covering~$\id_X$)
such that the restriction of~$\vph$ to each fiber is anti-complex linear
and $\vph\!\circ\!\vph\!=\!\id_V$.
A \sf{real bundle pair} $(V,\vph)\!\lra\!(X,\phi)$   
consists of a complex vector bundle $V\!\lra\!X$ and 
a conjugation~$\vph$ on $V$ lifting~$\phi$.
For example, 
$$(X\!\times\!\C^n,\phi\!\times\!\fc)\lra(X,\phi),$$
where $\fc\!:\C^n\!\lra\!\C^n$ is the standard conjugation on~$\C^n$,
is a real bundle pair.
If $X$ is a smooth manifold with a smooth involution~$\phi$,
then $(TX,\tnd\phi)$ is also a real bundle pair over~$(X,\phi)$.
For any real bundle pair $(V,\vph)$ over~$(X,\phi)$, the fixed locus~$V^{\vph}$
of~$\vph$ is a real vector bundle over~$X^{\phi}$.
We denote~by
$$\La_{\C}^{\top}(V,\vph)=\big(\La_{\C}^{\top}V,\La_{\C}^{\top}\vph\big)$$
the top exterior power of $V$ over $\C$ with the induced conjugation.
Direct sums, duals, and tensor products over~$\C$ of real bundle pairs over~$(X,\phi)$
are again real bundle pairs over~$(X,\phi)$.

\begin{dfn}[{\cite[Definition~5.1]{RealGWsI}}]\label{realorient_dfn4}
Let $(X,\phi)$ be a topological space with an involution and 
$(V,\vph)$ be a real bundle pair over~$(X,\phi)$.
A \sf{real orientation} on~$(V,\vph)$ consists~of 
\begin{enumerate}[label=(RO\arabic*),leftmargin=*]

\item\label{LBP_it2} a rank~1 real bundle pair $(L,\wt\phi)$ over $(X,\phi)$ such that 
\BE{realorient_e4}
w_2(V^{\vph})=w_1(L^{\wt\phi})^2 \qquad\hbox{and}\qquad
\La_{\C}^{\top}(V,\vph)\approx(L,\wt\phi)^{\otimes 2},\EE

\item\label{isom_it2} a homotopy class of isomorphisms of real bundle pairs in~\eref{realorient_e4}, and

\item\label{spin_it2} a spin structure~on the real vector bundle
$V^{\vph}\!\oplus\!2(L^*)^{\wt\phi^*}$ over~$X^{\phi}$
compatible with the orientation induced by~\ref{isom_it2}.\\ 

\end{enumerate}
\end{dfn}

\noindent
An isomorphism in~\eref{realorient_e4} restricts to an isomorphism 
$\La_{\R}^{\top}V^{\vph}\!\approx\!(L^{\wt\phi})^{\otimes2}$
of real line bundles over~$X^{\phi}$.
Since the vector bundles $(L^{\wt\phi})^{\otimes2}$ and $2(L^*)^{\wt\phi^*}$ are canonically oriented, 
\ref{isom_it2} determines orientations on $V^{\vph}$ and $V^{\vph}\!\oplus\! 2(L^*)^{\wt\phi^*}$.
By the first assumption in~\eref{realorient_e4}, the real vector bundle
$V^{\vph}\!\oplus\!2(L^*)^{\wt\phi^*}$ over~$X^{\phi}$ admits a spin structure.

\begin{prp}[{\cite[Proposition~7.3]{RBP}}]\label{canonisom_prp}
Let $(V,\vph)$ be a rank~$n$ real bundle pair over
a (possibly nodal) symmetric surface $(\Si,\si)$. 
A real orientation on $(V,\vph)$ determines a homotopy class of isomorphisms
\BE{canonisom_e} \Psi\!: \big(V\!\oplus\!2L^*,\vph\!\oplus\!2\wt\phi^*\big)
\approx\big(\Si\!\times\!\C^{n+2},\si\!\times\!\fc\big)\EE
of real bundle pairs over $(\Si,\si)$.
\end{prp}

\noindent
The existence of the isomorphisms~\eref{canonisom_e} over smooth symmetric surfaces
is implied by the classification of real bundle pairs over smooth symmetric surfaces in
\cite[Propositions~4.1,4.2]{BHH}. 
Proposition~\ref{canonisom_prp}, which is inspired by the direct sign computations 
in~\cite{XCapsSigns}, specifies topological data determining a homotopy class 
of such isomorphisms. 
This proposition provides the topological foundations for the approach of~\cite{RealGWsI}
to Problems~\ref{orient_it} and~\ref{bnd_it} on page~\pageref{orient_it}.
For the purposes of~\cite{RealGWsI}, it is sufficient to establish 
Proposition~\ref{canonisom_prp} for smooth and one-nodal symmetric surfaces;
these cases are \cite[Propositions~5.2,6.2]{RealGWsI}.
The case for symmetric surfaces with one pair of conjugate nodes is
\cite[Lemma~4.4]{RealGWsII}.
The principles behind the proofs of these special cases in~\cite{RealGWsI,RealGWsII}
are leveraged in~\cite{RBP} to give an alternative proof of \cite[Propositions~4.1,4.2]{BHH},
extend them to nodal symmetric surfaces, classify the homotopy classes of automorphisms
of real bundle pairs, and
establish the full statement of Proposition~\ref{canonisom_prp} as a corollary
of these results.\\

\noindent
If $(L,\wt\phi)$ is a rank~1 real bundle pair over a symmetric surface~$(\Si,\si)$
such that \hbox{$L^{\wt\phi}\!\lra\!\Si^{\si}$} is orientable, then the real vector bundle
$(L^{\wt\phi})^{\otimes2}\!\oplus\!2(L^*)^{\wt\phi^*}$ has a canonical spin structure.
Otherwise, such a spin structure can be canonically  fixed by a choice of orientation 
of each loop in~$\Si^{\si}$ over which $L^{\wt\phi}$ is not orientable
and it depends on this choice.
Combined with this observation, Proposition~\ref{canonisom_prp} yields the following
conclusion.

\begin{crl}[{\cite[Corollary~5.6]{RealGWsI}}]\label{canonisom_crl}
Let  $(L,\wt\phi)$ be a rank~1 real bundle pair over a symmetric surface~$(\Si,\si)$. 
If $L^{\wt\phi}\!\lra\!\Si^{\si}$ is orientable, there exists a canonical homotopy 
class of isomorphisms 
\BE{canonisomCRL_e}\big(L^{\otimes2}\!\oplus\!2L^*,\wt\phi^{\otimes2}\!\oplus\!2\wt\phi^*\big)
\approx\big(\Si\!\times\!\C^3,\si\!\times\!\fc\big)\EE
of real bundle pairs over $(\Si,\si)$.
In general, the canonical homotopy class of isomorphisms~\eref{canonisomCRL_e}
is determined by the choice of orientation for each loop in~$\Si^{\si}$
over which $L^{\wt\phi}$ is not orientable.
\end{crl}

\noindent
Our notion of real orientation on~$(X,\om,\phi)$ can be viewed 
as the real arbitrary-genus analogue of 
the notions of spin structure and relative spin structure of \cite[Definition~8.1.2]{FOOO}
in the  genus~0 open GW-theory.
These structures induce orientations on determinants of generalized Cauchy-Riemann operators
in the open setting and orient moduli spaces of $J$-holomorphic disks.
In some cases, they can be used to orient  moduli spaces of real $J$-holomorphic maps
from $\P^1$ with the standard orientation-reversing involution~$\tau_2$.
In~\cite{RealGWsI}, we show that a real orientation can be used to orient 
compactified moduli spaces of real $J$-holomorphic maps in arbitrary genera
whenever the ``complex" dimension of the target~$X$ is odd.

\section{Real Gromov-Witten theory}
\label{RealGW_sec}

\noindent
A \sf{real orientation on a  real symplectic manifold
$(X,\om,\phi)$} is a real orientation on the real bundle pair $(TX,\tnd\phi)$.
We call a real symplectic manifold $(X,\om,\phi)$ \sf{real-orientable} if 
it admits a real orientation.
The examples include $\P^{2n-1}$, $X_5$, many other projective complete intersections,
and simply-connected real symplectic Calabi-Yau and 
real Kahler Calabi-Yau manifolds with spin fixed locus;
see \cite[Propositions~1.2,1.4]{RealGWsIII}.

\begin{thm}[{\cite[Theorem~1.3]{RealGWsI}}]\label{orient_thm}
Let $(X,\om,\phi)$ be a real-orientable $2n$-manifold, 
$g,l\!\in\!\mathbb{Z}^{\ge0}$, 
\hbox{$B\!\in\!H_2(X;\Z)$}, and $J\!\in\!\cJ_{\om}^{\phi}$.  
\begin{enumerate}[label=(\arabic*),leftmargin=*]

\item\label{odd_it} If $n\!\not\in\!2\Z$, a real orientation on $(X,\om,\phi)$ orients 
$\ov\fM_{g,l}(X,B;J)^{\phi}$.
\item If $n\!\in\!2\Z$ and $g\!+\!l\!\ge\!2$,  
a  real orientation on $(X,\om,\phi)$ orients 
the real line~bundle
$$\La_{\R}^{\top}\big(T\ov\fM_{g,l}(X,B;J)^{\phi}\big)
\otimes\ff^*\La_\R^{\top}\big(T\R\ov\cM_{g,l}\big)
\lra \ov\fM_{g,l}(X,B;J)^{\phi}.$$

\end{enumerate}
\end{thm}

\noindent
Just as happens in   complex GW-theory,
an orientation on $\ov\fM_{g,l}(X,B;J)^{\phi}$ determined by some topological data on~$(X,\om,\phi)$
gives rise to a virtual class for this moduli space.
For each $i\!=\!1,\ldots,l$, let 
$$\ev_i\!: \ov\fM_{g,l}(X,B;J)^{\phi}\lra X, \qquad
\big[u,(z_1^+,z_1^-),\ldots,(z_l^+,z_l^-)\big]\lra u(z_i^+),$$
be the evaluation at the first point in the $i$-th pair of conjugate points.
For $\mu_1,\ldots,\mu_l\!\in\!H^*(X)$, the numbers
$$\blr{\mu_1,\ldots,\mu_l}\equiv \int_{[\ov\fM_{g,l}(X,B;J)^{\phi}]}
\ev_1^*\mu_1\ldots\ev_l^*\mu_l\in\Q$$
are virtual counts of real $J$-holomorphic curves in~$X$ passing through generic cycle
representatives for the Poincare duals of $\mu_1,\ldots,\mu_l$,
i.e.~\sf{real GW-invariants} of $(X,\om,\phi)$ with conjugate pairs of insertions.
They are independent of the choices of cycles representatives and of~$J$.\\

\noindent
As in  the complex GW-theory, it is convenient to consider moduli spaces of 
$J$-holomorphic maps from disconnected domains.
The topological components of a disconnected nodal symmetric surface~$(\Si,\si)$ 
split into those preserved by~$\si$ and into pairs of components interchanged by~$\si$;
the latter are called \sf{$g_0$-doublets} in \cite[Section~1.3]{RealGWsII},
where $g_0$ is the arithmetic genus of either topological component.
A real orientation on $(X,\om,\phi)$ with $n\!\not\in\!2\Z$ also determines
an orientation on the moduli spaces of real $J$-holomorphic maps from doublets.
Since any such map is determined by its restriction to either topological component
of the domain, 
the moduli spaces of real $J$-holomorphic maps from $g_0$-doublets with ordered components
can also be oriented from the standard complex orientation of 
the moduli space of  $J$-holomorphic genus~$g_0$ maps.
The two orientations differ by $(-1)^{g_0+1+l_-}$, where 
$l_-$ is the number of second points in each pair carried by the preferred component
of the doublet; see \cite[Theorem~1.4]{RealGWsII}.\\

\noindent
There are also alternative ways of orienting the uncompactified moduli spaces
$\fM_{0,l}(X,B;J)^{\phi,\si}$ for the two topological types of involutions on~$\P^1$
under the assumptions of Theorem~\ref{orient_thm}.
For the standard involution $\tau\!=\!\tau_2$, 
the orientation of Theorem~\ref{orient_thm} on $\fM_{0,l}(X,B;J)^{\phi,\tau}$
and the orientation induced as in \cite[Section~8.1]{FOOO} by the associated spin
or relative structure are the same up to a topologically determined sign;
see \cite[Theorem~1.5]{RealGWsII}.
For the fixed-point-free involution~$\eta$, 
the orientation of Theorem~\ref{orient_thm} on $\fM_{0,l}(X,B;J)^{\phi,\eta}$
and the orientation induced as in \cite[Section~2.1]{Teh} by the real square root
in~\eref{realorient_e4} are the same.
These comparisons extend to higher genus {\it once} $\fM_{g,l}(X,B;J)^{\phi,\tau}$
is known to be orientable; see \cite[Corollary~3.8]{RealGWsII}.\\

\noindent
There is a natural immersion
\BE{iomap_e}
\big\{\big[u,(z_1^+,z_1^-),\ldots,(z_{l+2}^+,z_{l+2}^-)\big]\!\in\!
\ov\fM_{g,l+2}^{\bu}(X,B;J)^{\phi}\!:\,
u(z_{l+1}^+)\!=\!u(z_{l+2}^+)\big\}\lra \ov\fM_{g+2,l}^{\bu}(X,B;J)^{\phi}\EE
from a subspace of the moduli space of maps from disconnected domains;
it identifies the marked points $z_{l+1}^{\pm}$ and $z_{l+2}^{\pm}$
to form a pair of conjugate nodes.
The analogue of this immersion in complex GW-theory identifies the last pair of points
to form a nodal surface of genus one higher.
The domain, target, and the normal bundle of this immersion in complex GW-theory
have natural complex orientations; this immersion is orientation-preserving
with respect to these orientations.
Under the conditions of Theorem~\ref{orient_thm}, a real orientation on~$(X,\om,\phi)$
induces orientations on the domain and target of the immersion~\eref{iomap_e};
the normal bundle of this immersion has a natural complex orientation.
By \cite[Theorem~1.2]{RealGWsII}, this immersion is orientation-{\it reversing}
with respect to these orientations.\\

\noindent
The forgetful morphism~\eref{ffdfn_e} and the comparisons of orientations 
in the last three paragraphs are essential parts of a fully fledged GW-theory 
needed for explicit computations of these invariants
and algebraic interpretations in the spirit of cohomological field theory.
In particular, they are used to describe equivariant localization data that computes
the arbitrary-genus real GW-invariants  of~$\P^{2n-1}$ with conjugate pairs of insertions
in \cite[Section~4.2]{RealGWsIII} and relate real GW- and enumerative invariants
 in \cite[Theorem~1.1]{NZ}.\\

\noindent
In some cases, we are able to define invariant counts of $J$-holomorphic curves
with real insertions as well.
However, similarly to the situation with \cite{Wel4,Wel6,Cho,Sol},
the relevant moduli spaces are not orientable in these cases and 
there is no fully fledged GW-theory involving these invariants.

\begin{thm}[{\cite[Theorem~1.5]{RealGWsI}}]\label{dim3_thm} 
Let $(X,\om,\phi)$ be a compact real-orientable 6-manifold 
such that $\lr{c_1(X),B}\!\in\!4\Z$ for all $B\!\in\!H_2(X;\Z)$
with $\phi_*B\!=\!-B$.
For all $B\!\in\!H_2(X;\Z)$ and $k,l\!\in\!\Z^{\ge0}$,
a real orientation on~$(X,\om,\phi)$ determines a  count 
$\lr{\pt^l;\pt^k}_{1,B}^{\phi}\!\in\!\Q$
of real $J$-holomorphic genus~1 degree~$B$ curves passing through generic collections
of $k$~real points and of $l$~pairs of conjugate points in~$X$.
This count is  independent of generic choices of the points and~$J\in\!\cJ_{\om}^{\phi}$.
\end{thm}


 \section{Orienting Fredholm determinants}
\label{MainPf_sec}

\noindent
A \textsf{real Cauchy-Riemann} (or \sf{CR-}) \sf{operator} on a real bundle pair $(V,\vph)$ 
over a symmetric Riemann surface $(\Si,\si,\fj)$ is a linear map of the~form
\begin{equation*}\begin{split}
D_V=\dbar_V\!+\!A\!: \Ga(\Si;V)^{\vph}
\equiv&\big\{\xi\!\in\!\Ga(\Si;V)\!:\,\xi\!\circ\!\si\!=\!\vph\!\circ\!\xi\big\}\\
&\hspace{.1in}\lra
\Ga^{0,1}(\Si;V)^{\vph}\equiv
\big\{\ze\!\in\!\Ga(\Si;(T^*\Si,\fj)^{0,1}\!\otimes_{\C}\!V)\!:\,
\ze\!\circ\!\tnd\si=\vph\!\circ\!\ze\big\},
\end{split}\end{equation*}
where $\dbar_V$ is the holomorphic $\dbar$-operator for some holomorphic structure in~$V$ and  
$$A\in\Ga\big(\Si;\Hom_{\R}(V,(T^*\Si,\fj)^{0,1}\!\otimes_{\C}\!V) \big)^{\vph}$$ 
is a zeroth-order deformation term. 
Let $\dbar_{\Si;\C}$ denote the real CR-operator on 
the trivial rank~1 real bundle
\hbox{$(\Si\!\times\!\C,\si\!\times\fc)$}
with the standard holomorphic structure and $A\!=\!0$.\\

\noindent
Any real CR-operator~$D_V$ on a real bundle pair $(V,\vph)$  
over a symmetric Riemann surface $(\Si,\si,\fj)$ is Fredholm in the appropriate completions.
We denote~by
$$\det D_V\equiv\La_{\R}^{\top}(\ker D_V) \otimes \big(\La^{\top}_{\R}(\text{cok}\,D_V)\big)^*$$
its \sf{determinant line}.
Since the space of real CR-operators on $(V,\vph)$ is contractible,
an orientation on $\det D_V$ for one such operator determines an orientation for
 all real CR-operators on~$(V,\vph)$.
Thus, an exact sequence
$$0\lra (V_1,\vph_1)\lra (V,\vph)\lra (V_2,\vph_2)\lra 0$$
of real bundle pairs over $(\Si,\si)$ determines a homotopy class of isomorphisms
\BE{detisom_e} \det D_V \approx \big(\!\det D_{V_1}\big)\otimes \big(\!\det D_{V_2}\big)\EE
between the determinants of any real CR-operators on these real bundle pairs.
Via these isomorphisms, orientations on any two of the determinants in~\eref{detisom_e}
determine an orientation on the third.
Furthermore, the line $(\det D_V)^{\otimes2}$ is canonically oriented for any real CR-operator~$D_V$.
By Proposition~\ref{canonisom_prp},  a real orientation on a rank~$n$ real bundle 
pair~$(V,\vph)$ over $(\Si,\si)$ thus determines an orientation on the~line
\BE{fDdfn_e}
\rdet\,D_V\equiv \big(\!\det D_V\big)\otimes\big(\!\det\dbar_{\Si;\C}\big)^{\otimes n}\EE
for every real CR-operator~$D_V$  on the real bundle pair $(V,\vph)$ over~$(\Si,\si)$.
We call $\rdet\,D_V$ \sf{the relative determinant of~$D_V$},
since an orientation on $\rdet\,D_V$ determines a correspondence between 
the orientations on $\det D_V$ and on the determinant of $\det n\dbar_{\Si;\C}$
of the standard real CR-operator on 
the trivial rank~$n$ real bundle
\hbox{$(\Si\!\times\!\C^n,\si\!\times\fc)$} over~$(\Si,\si)$.\\

\noindent
For each element $[u]$ of $\ov\fM_{g,l}(X,B;J)^{\phi}$, the linearization
$$D_u\!: \Ga(\Si;u^*TX)^{u^*\tnd\phi}\lra \Ga^{0,1}(\Si;u^*TX)^{u^*\tnd\phi}$$ 
of the real $\dbar_J$-operator at~$[u]$ is a real CR-operator.
If $g\!+\!l\!\ge\!2$, the forgetful morphism~\eref{ffdfn_e} induces a canonical isomorphism
\BE{thm_maps_e3}\La_{\R}^{\top}\big(T_{[u]}\ov\fM_{g,l}(X,B;J)^{\phi}\big)
\approx \big(\!\det D_u\big)\otimes \La_{\R}^{\top}(T_{\ff([u])}\R\ov\cM_{g,l}).\EE
Theorem~\ref{orient_thm}\ref{odd_it} may then appear to be about systematically 
orienting  each factor on the right-hand side of~\eref{thm_maps_e3}.
However, the moduli space $\R\ov\cM_{g,l}$ is not orientable if $g\!\in\!\Z^+$;
a systematic orientation on the left-hand side of~\eref{thm_maps_e3} exists 
if and only if the family $\det D_u$ is not orientable in the same manner.\\

\noindent
A real orientation on~$(X,\om,\phi)$ pulls back to a real orientation
on $u^*(TX,\tnd\phi)$ and thus systematically determines an orientation on the relative 
determinant
\BE{fDdfn_e2}
\rdet\,D_u\equiv \big(\!\det D_u\big)\otimes\big(\!\det\dbar_{\Si;\C}\big)^{\otimes n}\EE
of~$D_u$ for each element $[u]$ of $\ov\fM_{g,l}(X,B;J)^{\phi}$.
It is immediate that this orientation varies continuously over each stratum of 
 $\ov\fM_{g,l}(X,B;J)^{\phi}$ and 
straightforward to show that it varies continuously across the strata as well;
see \cite[Corollary~6.7]{RealGWsI} for the crucial extension from 
the main stratum across the codimension~1 strata.
In this light, the last factor of~\eref{fDdfn_e2} describes the orientability of~$\det\,D_u$.\\

\noindent
The orientability of $\R\ov\cM_{g,l}$ is described
by Proposition~\ref{DM_prp} below.
We denote~by 
$$\det\dbar_{\C}\lra \R\ov\cM_{g,l}$$
the line bundle with fiber $\det\dbar_{\Si;\C}$ over~$[\Si]$.

\begin{prp}[{\cite[Propositions~5.9,6.1]{RealGWsI}}]\label{DM_prp} 
Let $g,l\in \Z^{\ge0}$ be such that $g\!+\!l\!\ge\!2$.
The restriction of the line bundle 
\BE{CidentDM_e}\La_{\R}^{\top}\big(T\R\ov\cM_{g,l}^{\si}\big)
\otimes\big(\!\det\dbar_{\C}\big)\lra \R\ov\cM_{g,l}\EE
to the main stratum $\R\cM_{g,l}$ consisting of smooth symmetric surfaces is
canonically oriented.
This canonical orientation extends over $\R\ov\cM_{g,l}$ after 
it is reversed
over every topological component $\cM_{g,l}^{\si}\!\subset\!\R\cM_{g,l}$
of smooth Riemann surfaces~$(\Si,\si,\fj)$ with $g\!-\!|\pi_0(\Si^{\si})|\!\in\!2\Z$.
\end{prp}

\noindent
The first statement of this proposition is obtained by combining the Kodaira-Spencer (\sf{KS})
and Serre Duality (\sf{SD}) isomorphisms with the first claim of Corollary~\ref{canonisom_crl} 
for $(L,\wt\phi)\!=\!(T^*\Si,(\tnd\si)^*)$.
Each of these three isomorphisms induces an orientation on the restriction of
a real line bundle over $\R\ov\cM_{g,l}$ to $\R\cM_{g,l}$;
the tensor product of these lines bundles yields the desired orientation.
This first statement, which concerns only smooth symmetric surfaces,
 is also obtained in~\cite{Remi}, but
based on the highly technical analysis of the action of automorphisms of real bundle pairs
on the determinants of real CR-operators in~\cite{Remi0} instead of 
the topological claim of Corollary~\ref{canonisom_crl}.\\

\noindent
Corollary~\ref{canonisom_crl} is further used to establish the second claim of
Proposition~\ref{DM_prp} and thus Theorem~\ref{orient_thm}.
The orientation induced by the KS~isomorphism flips across 
the codimension~1 boundary strata of  $\R\ov\cM_{g,l}$,
while the orientation induced by the SD~isomorphism extends across these strata.
The orientation induced by  the first claim of Corollary~\ref{canonisom_crl} 
extends across the codimension~1 boundary strata of types~(E) and~(H1) on page~\pageref{E_it}
and flips across the codimension~1 boundary strata of types~(H2) and~(H3).
This is an artifact of the topology of the continuous extension~$(\wh\cT,\wh\vph)$
of $(T^*\Si,(\tnd\si)^*)$ to a one-nodal symmetric surface~$(\Si,\si)$.
If $(\Si,\si)$ is of type~(E) or~(H1), then $\wh\cT^{\wh\vph}$ is orientable
and the first claim of Corollary~\ref{canonisom_crl} applies;
an isomorphism~\eref{canonisomCRL_e} in the canonical homotopy class extends  
to isomorphisms in the canonical homotopy classes for the nearby symmetric surfaces.
If $(\Si,\si)$ is of type~(H2) and~(H3), then $\wh\cT^{\wh\vph}$ not orientable
and there is no choice
of orientation of the loops in~$\Si^{\si}$ which extends to both possible smoothing directions;
the second claim of Corollary~\ref{canonisom_crl} applies in this case
and is responsible for the change in the orientation.\\

\begin{table}
\begin{center}
\begin{small}
\renewcommand\arraystretch{1.8} 
\begin{tabular}{||c|c|c||}
\hline\hline
& (E)/(H1)& (H2)/(H3)\\
\hline
orientation on~\eref{fDdfn_e2} with $(V,\vph)\!=\!u^*(TX,\tnd\phi)$&  $+$ & $+$\\
\hline
orientation induced by KS isomorphism& $-$ & $-$\\
\hline
orientation induced by SD isomorphism& $+$ & $+$\\
\hline
orientation on~\eref{fDdfn_e2} with $(V,\vph)\!=\!(T^*\Si,(\tnd\si)^*)^{\otimes2}$ & $+$ & $-$\\
\hline\hline
parity of $|\pi_0(\Si^{\si})|$ &  $-$ & $+$\\
\hline\hline
\end{tabular}
\end{small}
\end{center}
\caption*{\it The extendability of the canonical orientations factoring 
into~\eref{twistedT_e}
and of the parity of 
the number of components of~$\Si^{\si}$ across the codimension~1 strata:
$+$ extends, $-$ flips.}
\end{table}

\noindent
By the previous paragraph,  
the canonical orientation on the restriction of the line bundle~\eref{CidentDM_e}
to $\R\cM_{g,l}$ flips across the codimension~1 boundary strata of types~(E) and~(H1)
and extends across the codimension~1 boundary strata of types~(H2) and~(H3).
The parity of $|\pi_0(\Si^{\si})|$ behaves in the same way.
These two statements together establish the second claim of Proposition~\ref{DM_prp}.\\

\noindent
By \eref{thm_maps_e3}, \eref{fDdfn_e2}, and Proposition~\ref{DM_prp},
a real orientation on~$(X,\om,\phi)$ determines an orientation on 
\BE{twistedT_e}\La_{\R}^{\top}\big(T_{[u]}\ov\fM_{g,l}(X,B;J)^{\phi}\big)
\otimes \big(\!\det\dbar_{\Si;\C}\big)^{\otimes(n+1)}
\approx  \big(\rdet\,D_u\big)\otimes
\La_{\R}^{\top}\big(T_{[\Si]}\cM_{g,l}^{\si}\big)
\otimes\big(\!\det\dbar_{\Si;\C}\big)\EE
which depends continuously on $[u]\!\in\!\ov\fM_{g,l}(X,B;J)^{\phi}$.
If $n\!\not\in\!2\Z$, an orientation on~\eref{twistedT_e}
is equivalent to an orientation on $T_{[u]}\ov\fM_{g,l}(X,B;J)^{\phi}$.
This establishes Theorem~\ref{orient_thm} whenever $g\!+\!l\!=\!2$;
the three exceptional cases are then deduced by adding extra conjugate pairs of marked points.

\section{Real enumerative geometry}
\label{RealEnum_sec}

\noindent
As in the complex case, the curve-counting invariants arising from Theorem~\ref{orient_thm}
are generally rational numbers.
For specific real almost Kahler manifolds $(X,\om,\phi,J)$,
they can be converted into signed counts of genus~$g$ degree~$B$ real 
$J$-holomorphic curves passing through specified conjugate pairs of constraints 
and thus provide lower bounds in real enumerative geometry.
If $J$ is a sufficiently regular almost complex structure on $(X,\om,\phi)$,
then the two curve counts are the same in genus~0
(this is also the case in the complex setting).
In the real (but not complex) setting, this equality extends to 
the genus~1 curve counts in 6-folds.

\begin{prp}[{\cite[Theorem~1.5]{RealGWsIII}}]\label{g1EG_thm}
Let $(X,\om,\phi)$ be a compact real-orientable 6-fold 
and $J\!\in\!\cJ_{\om}^{\phi}$ be a generic 
almost complex structure on $(X,\om,\phi)$.
The genus~1 real GW-invariants of $(X,\om,\phi)$ arising from
Theorems~\ref{orient_thm} and~\ref{dim3_thm} are then equal to 
the corresponding signed counts of real $J$-holomorphic curves
and thus provide lower bounds for the number of real genus~1 irreducible curves in~$(X,J,\phi)$. 
\end{prp}

\noindent
This statement extends to higher genus as follows.
The real GW-invariants arising from Theorem~\ref{orient_thm} induce homomorphisms
$$\GW_{g,B}^{X,\phi}\!: H^*(X;\Z)^{\oplus l}\lra \Q$$
obtained by pulling back cohomology classes on~$X$ by the evaluation maps for  
the first marked points in the conjugate pairs.
For $g,h\!\in\!\Z^{\ge 0}$ and $B\!\in\!H_2(X;\Z)$, define $\wt{C}_{h,B}^X(g)\!\in\!\Q$ by
$$ \sum_{g=0}^{\i}\wt{C}_{h,B}^X(g)t^{2g}
=\bigg(\frac{\sinh(t/2)}{t/2}\bigg)^{\!h-1+\lr{c_1(X,\om),B}/2}\,.$$
Since $\wt{C}_{h,B}^X(0)\!=\!1$, we can define homomorphisms
\BE{FanoGV_e} 
\E_{h,B}^{X,\phi}\!: H^*(X;\Z)^{\oplus l}\lra\Q ~~~\forall\,h\!\in\!\Z^{\ge 0}
\quad\hbox{by}\quad 
\GW_{g,B}^{X,\phi}=
\sum_{\begin{subarray}{c}0\le h\le g\\ g-h\in 2\Z\end{subarray}}\!\!\!\!
\wt{C}_{h,B}^X\big(\tfrac{g-h}{2}\big) \E_{h,B}^{X,\phi}
~~~\forall\,g\!\in\!\Z^{\ge0}.\EE

\begin{thm}[{\cite[Theorem~1.1]{NZ}}]\label{GWvsEnum_thm}
Suppose $(X,\om,\phi)$ is a compact real-orientable symplectic \hbox{6-fold} 
with a choice of real orientation, $g,l\!\in\!\Z^{\ge 0}$, and $B\!\in\!H_2(X;\Z)$.
If $\lr{c_1(\om,X),B}\!>\!0$, then the homomorphisms $\E_{h,B}^{X,\phi}$ take values in~$\Z$.
If in addition $J$ is a generic almost complex structure on $(X,\om,\phi)$ and 
$\mu_1,\ldots,\mu_l\!\in\!H^*(X;\Q)$ are such~that
$$\sum_{i=1}^l\dim_\R\mu_i=\lr{c_1(X,\om),J}+2l,$$
then  $\E_{h,B}^{X,\phi}(\mu_1,\ldots,\mu_l)$ is the number of real irreducible $J$-holomorphic 
genus~$h$  degree~$B$ curves 
passing through a generic collection of cycles representing  $\mu_1,\ldots,\mu_l$
counted with~sign.
\end{thm}

\noindent
This is the real analogue of \cite[Theorem~1.5]{FanoGV}
which established the ``Fano" case of the Gopakumar-Vafa prediction 
of \cite[Conjecture~2(i)]{P2} in a stronger~form.
In the case of~$\P^3$, there is a natural identification of the second $\Z$-homology
with~$\Z$.
By~\eref{FanoGV_e} and \cite[Theorem~1.6]{RealGWsIII}, $\E_{h,d}^{X,\phi}\!=\!0$
whenever $d\!-\!h\!\in\!2\Z$.
The standard complex structure of $\P^3$ is ``generic" for the purposes 
of Theorem~\ref{GWvsEnum_thm} if $d\!\ge\!2h\!-\!1$;
it is ``generic" for some lower values of~$d$ as well.
The equivariant localization data of \cite[Section~4.2]{RealGWsIII}
is used in \cite{NZapp} to compute the genus~$g$ degree~$d$ real GW-invariants of~$\P^3$
with $d$~conjugate pairs of point insertions for $g\!\le\!5$ and $d\!\le\!8$.
The real enumerative invariants obtained from these numbers via~\eref{FanoGV_e}
and shown in \cite[Table~2]{NZ} are consistent with the complex enumerative invariants
and thus with the Castelnuovo bounds.\\

\noindent
The number of real curves passing through specified constraints generally depends
on the constraints themselves and not just on the (co)homology classes they represent;
only a properly signed count of such curves can be invariant in general.
It is bounded above by the analogous complex count and has the same parity 
as the latter.
There are known examples when the upper bounds provided by the complex counts are sharp,
i.e.~there are nonempty open subsets of the spaces of admissible constraints 
achieving these bounds.
The real GW-invariants of \cite{Wel4,Wel6,Ge2,Teh} in genus~0 and of~\cite{RealGWsI}
in positive genera lead to lower bounds for counts of real curves in certain cases;
 examples when these bounds are sharp are obtained in~\cite{Kollar}.
It remains an open question when these upper and lower are sharp in general.\\

\vspace{.5in}

\noindent
{\it  Institut de Math\'ematiques de Jussieu - Paris Rive Gauche,
Universit\'e Pierre et Marie Curie, 
4~Place Jussieu,
75252 Paris Cedex 5,
France\\
penka.georgieva@imj-prg.fr}\\

\noindent
{\it Department of Mathematics, Stony Brook University, Stony Brook, NY 11794\\
azinger@math.stonybrook.edu}\\


\end{document}